\documentclass{gtart_h}

\def\ifplaintex{\expandafter\ifx\csname documentclass\endcsname\relax}

\def\gtp{{\mathsurround=0pt\it $\cal G\mskip-2mu$eometry \&\ 
$\cal T\!\!$opology $\cal P\!$ublications}}  

\def\recd{{\small Received:\qua\receiveddate\ifx\reviseddate\relax
\else\qquad Revised:\qua\reviseddate\fi\par}} 

\def\lognumber#1{\def\thelognumber{#1}}
\def\volumenumber#1{\def\thevolumenumber{#1}}
\def\volumeyear#1{\def\thevolumeyear{#1}}
\def\papernumber#1{\def\thepapernumber{#1}}
\def\pagenumbers#1#2{\def\startpage{#1}\def\finishpage{#2}}
\def\published#1{\def\publishdate{#1}}

\def\received#1{\def\receiveddate{#1}}

\def\accepted#1{\def\accepteddate{#1}}
\def\asciititle#1{\def\theasciititle{#1}}

\long\def\asciiabstract#1{\long\def\theasciiabstract{#1}}
\def\asciikeywords#1{\def\theasciikeywords{#1}}

\let\\\par\let\thelognumber\relax\let\thevolumenumber\relax
\let\thepapernumber\relax\let\thevolumeyear\relax\let\startpage\relax
\let\finishpage\relax\let\publishdate\relax\let\receiveddate\relax
\let\reviseddate\relax\let\accepteddate\relax\let\theasciititle\relax
\let\theasciiauthors\relax
\let\theasciiabstract\relax\let\theasciikeywords\relax

\let\theasciiemail\relax

\ifplaintex
\font\logobig=cmssbx10 scaled 3836
\font\logomed=cmssbx10 scaled 2557
\else
\font\logobig=cmssbx10 scaled 4200
\font\logomed=cmssbx10 scaled 2800
\fi

\long\def\makeagttitle{   
\count0=\startpage
\agt\hfill      
\hbox to 45truept{\vbox to 0pt{\vglue -13truept{\logomed A\kern -.37em{\logobig 
T}\kern -.38em G}\vss}\hss}
\break
{\small Volume \thevolumenumber\ (\thevolumeyear)
\startpage--\finishpage\nl
Published: \publishdate}

\vglue .25truein

{\parskip=0pt\leftskip 0pt plus
1fil\def\\{\par\smallskip}{\Large\bf\thetitle}\par\medskip} \vglue
0.05truein

{\parskip=0pt\leftskip 0pt plus 1fil\def\\{\par}{\sc\theauthors}
\par\medskip}%
 
\vglue 0.03truein 

{\small\leftskip 25truept\rightskip 25truept{\bf Abstract}\stdspace\theabstract

{\bf AMS Classification}\stdspace\theprimaryclass
\ifx\thesecondaryclass\relax\else; \thesecondaryclass\fi\par
{\bf Keywords}\stdspace \thekeywords\par}\vglue 7truept

}   

\ifplaintex
\font\phead=cmsl9 scaled 950
\font\pnum=cmbx10 scaled 913
\font\pfoot=cmsl9 scaled 950
\headline{\vbox to 0pt{\vskip -4.5mm\line{\small\phead\ifnum
\count0=\startpage ISSN 1472-2739 (on-line) 1472-2747 (printed)
\hfill {\pnum\folio}\else\ifodd\count0\def\\{ }%
\ifx\theshorttitle\relax\thetitle\else\theshorttitle\fi\hfill{\pnum\folio}
\else\def\\{ and }{\pnum\folio}\hfill\ifx\theshortauthors\relax\theauthors
\else\theshortauthors\fi\fi\fi}\vss}}
\footline{\vbox to 0pt{\vglue 0mm\line{\small\pfoot\ifnum\count0=\startpage
\copyright\ \gtp\hfill\else
\agt, Volume \thevolumenumber\ (\thevolumeyear)\hfill\fi}\vss}}
\else
\font=cmsl9 scaled 1050
\font\lnum=cmbx10 
\font=cmsl9 scaled 1050
\makeatletter
\def\@oddhead{{\small\lhead\ifnum\count0=\startpage ISSN 1472-2739 
(on-line) 1472-2747 (printed)\hfill {\lnum\number\count0}\else\ifodd\count0
\def\\{ }\ifx\theshorttitle\relax \thetitle \else\theshorttitle\fi\hfill
{\lnum\number\count0}\else\def\\{ and }{\lnum\number\count0}
\hfill\ifx\theshortauthors\relax 
\theauthors\else\theshortauthors\fi\fi\fi}}\def\@evenhead{\@oddhead}
\def\@oddfoot{\small\lfoot\ifnum\count0=\startpage\copyright\ \gtp\hfill\else
\agt, Volume \thevolumenumber\ (\thevolumeyear)\hfill\fi}
\def\@evenfoot{\@oddfoot}
\makeatother
\fi
\let\maketitlepage\makeagttitle
\let\makeshorttitle\maketitlepage
\let\maketitle\maketitlepage

\newwrite\gtoutfile
\long\gdef\makeheadfile{  
{\def\\{, }\def\s{ }
\immediate\openout\gtoutfile head.xxx
\immediate\write\gtoutfile{Proxy-for: \ifx\theasciiauthors\relax
\theauthors\else\theasciiauthors\fi\s<\ifx\theasciiemail\relax\theemail\else\theasciiemail\fi>}
\immediate\write\gtoutfile{\noexpand\\}
\immediate\write\gtoutfile{Authors: \ifx\theasciiauthors\relax
\theauthors\else\theasciiauthors\fi}
{\def\\{ }\immediate\write\gtoutfile{Title: \ifx\theasciititle\relax
\thetitle\else\theasciititle\fi}}
\immediate\write\gtoutfile{Subj-class: GT or SG, GR etc}
\immediate\write\gtoutfile{MSC-class: \theprimaryclass\ifx\thesecondaryclass\relax\else, \thesecondaryclass\fi}
\immediate\write\gtoutfile{Journal-ref: Algebr. Geom. Topol. \thevolumenumber\s
(\thevolumeyear) \startpage-\finishpage}
\immediate\write\gtoutfile{Comments: Published by Algebraic and
Geometric Topology at}
\immediate\write\gtoutfile{\s\s\s  http://www.maths.warwick.ac.uk/agt/AGTVol\thevolumenumber/agt-\thevolumenumber-\thepapernumber.abs.html}
\immediate\write\gtoutfile{\noexpand\\}
\immediate\write\gtoutfile{}
\ifx\theasciiabstract\relax
\immediate\write\gtoutfile{\theabstract}\else
\immediate\write\gtoutfile{\theasciiabstract}\fi
\immediate\write\gtoutfile{}
\immediate\write\gtoutfile{\noexpand\\}
\immediate\write\gtoutfile{}
\immediate\closeout\gtoutfile}}  

\def\maketitlepage{\makeagttitle\makeheadfile}
\let\makeshorttitle\maketitlepage
\let\maketitle\maketitlepage

\lognumber{66}
\volumenumber{5}
\volumeyear{2005}
\papernumber{66}
\pagenumbers{1637}{1653}
\received{16 September 2005} 
\accepted{21 November 2005}
\published{25 November 2005}

\usepackage{graphicx,amssymb,amsmath,psfrag,labelfig}

\psfrag {1}{\small1}
\psfrag {2}{\small2}
\psfrag {3}{\small3}
\psfrag {4}{\small4}
\psfrag {5}{\small5}
\psfrag {6}{\small6}
\psfrag {7}{\small7}
\psfrag {8}{\small8}
\psfrag {9}{\small9}

\def\figref#1{\hyperlink{#1anchor}{Figure~\ref*{#1}}}
\def\anchor#1{\noindent\hypertarget{#1anchor}{\smash{$\phantom{99}$}}\newline}

\newtheorem{proposition}{Proposition}
\newtheorem{theorem}[proposition]{Theorem}
\newtheorem{lemma}[proposition]{Lemma}
\newtheorem{corollary}[proposition]{Corollary}

\theoremstyle{definition}
\newtheorem{definition}[proposition]{Definition}
\newtheorem{step}{Step}

\def\v{\mathbf{v}}
\def\Kh{\mathit{Kh}}
\def\HKh{\mathit{HKh}}
\def\CKh{\mathit{CKh}}
\def\HKhshift{\mathit{HKh}_{\mathrm{sh}}}
\def\CKhshift{\mathit{CKh}_{\mathrm{sh}}}
\def\tb{tb}
\def\Res{\operatorname{Res}}
\def\mindeg{\operatorname{min-deg}}
\def\Z{\mathbb{Z}}
\def\Q{\mathbb{Q}}
\def\R{\mathbb{R}}
\def\s{\mathfrak{s}}
\def\notextit{}

\input xy
\xyoption{all}

\begin{document}

\title{A
Legendrian Thurston--Bennequin bound\\from Khovanov homology}
\asciititle{A Legendrian Thurston-Bennequin bound from Khovanov homology}
\author{Lenhard Ng}
\address{Department of Mathematics, Stanford University, Stanford,
CA 94305, USA}
\email{lng@math.stanford.edu}
\urladdr{http://alum.mit.edu/www/ng}

\begin{abstract}
We establish an upper bound for the Thurston--Bennequin number of a
Legendrian link using the Khovanov homology of the underlying
topological link. This bound is sharp in particular for all
alternating links, and knots with nine or fewer crossings.
\end{abstract}

\asciiabstract{%
We establish an upper bound for the Thurston-Bennequin number of a
Legendrian link using the Khovanov homology of the underlying
topological link. This bound is sharp in particular for all
alternating links, and knots with nine or fewer crossings.}


\asciikeywords{Legendrian link, Thurston-Bennequin number, Khovanov
homology, alternating link}

\maketitle

\section{Introduction}
\label{sec:intro}

The \textit{standard contact structure} on $\R^3$ is the two-plane
distribution given by the kernel of the one-form $dz-y\,dx$. A
\textit{Legendrian link} in standard contact $\R^3$ is a link which
is everywhere tangent to this contact structure. It is well-known
that any topological link type has a Legendrian representative.

To any oriented Legendrian link (or unoriented Legendrian knot),
there is an invariant called the \textit{Thurston--Bennequin
number}, abbreviated $\tb$, which measures the framing of the
contact plane field around the link. Given any Legendrian link, one
can construct a Legendrian link in the same topological class but
with $\tb$ less by $1$, by a construction known as stabilization. On
the other hand, it is not always possible to increase $\tb$ within a
link type. A classic result of Bennequin \cite{bib:Ben} states that
$\tb$ is bounded above by minus the Euler characteristic of a
Seifert surface for the link. The Bennequin bound on $\tb$ is one of
the fundamental results in three-dimensional contact topology, and
implies for example the existence of a contact structure on $\R^3$
which is homotopic but not isomorphic to the standard one.

Since Bennequin's result, there has been considerable interest in
computing or bounding $\overline{\tb}(K)$, the maximal
Thurston--Bennequin number for Legendrian links in the topological
link type $K$. Upper bounds on $\overline{\tb}$ in terms of other
knot invariants come from the following inequalities, where we abuse
notation and use $K$ to denote both a Legendrian link and its
underlying topological link:
\begin{itemize}
\item
$\tb(K) + |r(K)| \leq 2g(K)-1$, Bennequin's original result
\cite{bib:Ben}, where $g$ is the genus of $K$;
\item
$\tb(K) + |r(K)| \leq 2g_4(K) - 1$ \cite{bib:Rud2}, where $g_4$ is
the slice genus of $K$;
\item
$\tb(K) + |r(K)| \leq \mindeg_a H_K(a,z) - 1$ \cite{bib:FW,bib:Mo},
where $H_K$ is the HOMFLYPT polynomial and $\mindeg_a$ is the
minimum degree in the variable $a$;
\item
$\tb(K) + |r(K)| \leq 2\tau(K)-1$ \cite{bib:Pla}, where $\tau(K)$ is
the concordance invariant from knot Floer homology \cite{bib:OSz};
\item
$\tb(K) + |r(K)| \leq s(K)-1$ \cite{bib:Pla2,bib:Shu2}, where $s(K)$
is Rasmussen's $s$-invariant \cite{bib:Ras};
\item
$\tb(K) + |r(K)| \leq \hat{g}_{\text{min}}(K)$ \cite{bib:Wu}, where
$\hat{g}_{\text{min}}(K)$ is the minimal $k$ such that the
Khovanov--Rozansky cohomology $\mathit{HKR}^j_{k,l}(K)$
\cite{bib:KhR} is nonvanishing for some $j,l$;
\item
$\tb(K) \leq \mindeg_a F_K(a,z) - 1$ \cite{bib:Rud1}, where $F_K$ is
the Kauffman polynomial.
\end{itemize}

\noindent See \cite{bib:Fer} for further history (up to 2002). In
the first six bounds, $r(K)$ is the rotation number of Legendrian
$K$. Note that $r(K)$ is present in the first six bounds because
they in fact bound the self-linking number of a transverse link. Any
bound on $\tb+|r|$ from transverse knot theory is fundamentally
different from a direct bound on $\tb$ such as the Kauffman bound
above. For instance, the maximal value for $\tb+|r|$ for a left
handed trefoil is $-5$; no bound on $\tb+|r|$ can give the sharp
result $\overline{\tb}=-6$ in this case.

It is thus not surprising that the Kauffman bound tends to be more
effective in bounding $\overline{\tb}$ for Legendrian links than the
other bounds. The Kauffman bound is sharp (i.e., equality is
attained) for two-bridge links, and for all but two knots with
crossing number at most $9$ \cite{bib:NgTB}. It is, however, not
sharp for many negative torus knots, for which $\overline{\tb}$ has
been computed by Etnyre and Honda \cite{bib:EH} using
symplectic-topological techniques.

In this paper, we establish a new bound for $\overline{\tb}$ in
terms of ($\mathfrak{sl}(2)$) Khovanov homology \cite{bib:Kh}. This
bound is sharp for all alternating links, as well as all knots with
crossing number at most $9$, and all but at most two $10$-crossing
knots. In general, it seems to give the best overall currently known
bound on $\overline{\tb}$, although it is still not sharp for many
negative torus knots, and occasionally gives a worse bound than
Kauffman.

Recall that Khovanov homology associates to any oriented link $K$ a
bigraded abelian group $\HKh^{*,*}(K)$, where the first grading is
called the quantum (or Jones) grading and the second is the
homological grading. If we disregard torsion, we obtain a
two-variable Poincar\'e polynomial
\[
\Kh_K(q,t) = \sum \dim_{\Q}(\HKh^{i,j}(K)\otimes\Q) q^i t^j.
\]
Khovanov homology is designed so that its graded Euler
characteristic is the Jones polynomial $V_K$: $\Kh_K(q,-1) =
(q+q^{-1}) V_K(q^2)$.

\begin{theorem}[strong Khovanov bound]
For any link $K$,
\label{thm:strong}
\[
\overline{\tb}(K) \leq \min \{k\,|\, \bigoplus_{i-j=k} \HKh^{i,j}(K)
\neq 0\}.
\]
\end{theorem}

\begin{corollary}[weak Khovanov bound]
$\overline{\tb}(K) \leq \mindeg_q \Kh_K(q,t/q)$.
\label{cor:weak}
\end{corollary}

Theorem~\ref{thm:strong} and Corollary~\ref{cor:weak} often given
the same bound on $\overline{\tb}$, but there are instances in which
Theorem~\ref{thm:strong} is stronger; see
Section~\ref{sec:computations}. We may deduce from the Khovanov
bound another proof of the Legendrian corollary of the $s$-invariant
bound above:

\begin{corollary}
If $K$ is a knot, then $\overline{\tb}(K) \leq s(K)-1$.
\label{cor:s-invt}
\end{corollary}

\noindent Note that Corollary~\ref{cor:s-invt}, in turn, implies the
slice Bennequin bound $\overline{\tb}(K) \leq 2g_4(K)-1$; however,
it is generally nowhere as effective as the Khovanov bound in
bounding $\overline{\tb}$.

By using the Khovanov bound, we can calculate $\overline{tb}$ for
alternating links.

\begin{theorem}
The Khovanov bound (either weak or strong) is sharp for alternating
links. If $K$ is alternating and nonsplit, and $V_K(t),\sigma(K)$
are the Jones polynomial and signature of $K$ respectively, then
\label{thm:alternating}
\[
\overline{\tb}(K) = \mindeg_q V_K(q) + \sigma(K)/2 - 1.
\]
Here we use the convention that the signature of the right handed
trefoil is $+2$.
\end{theorem}

The reader may have noticed that the grading collapse in the strong
Khovanov bound is the same as the one which Seidel and Smith use in
their construction of a link invariant from Lagrangian intersection
Floer homology \cite{bib:SS}. Indeed, it is likely that one could
bound $\overline{\tb}(K)$ above by the minimum degree in which
Seidel--Smith's ``symplectic Khovanov cohomology''
$\HKh_{\mathit{symp}}^*(K)$ does not vanish. (Here $K$ might need to
be replaced by its mirror, depending on conventions.) This would use
unpublished work of Lipshitz and Manolescu constructing generators
for symplectic Khovanov cohomology in terms of a bridge diagram of a
link, generalizing \cite{bib:Man}. However, such an upper bound
would also follow from the strong Khovanov bound, if, as proposed in
\cite{bib:SS}, there is a spectral sequence from $\HKh$ to
$\HKh_{\mathit{symp}}$; in particular, Seidel--Smith's conjecture
that $\HKh^k_{\mathit{symp}}(K) \cong \bigoplus_{i-j=k}
\HKh^{i,j}(K)$ would imply that the two bounds are actually
identical.

We prove Theorem~\ref{thm:strong} in Section~\ref{sec:KhBound}. The
method of proof, which uses induction, seems to be somewhat
different from the proofs of previous Thurston--Bennequin bounds. In
particular, it yields a sufficient condition for a front to maximize
$\tb$, which can then be applied to calculate $\overline{\tb}$ for
alternating knots and nonsplit links. The construction of
maximal-$\tb$ fronts for alternating links requires some graph
theory which is the subject of Section~\ref{sec:alternating}. We
examine the Khovanov bound in some other illustrative examples in
Section~\ref{sec:computations}.

\rk{Acknowledgments}
I thank Robert Lipshitz and Ciprian Manolescu for useful
discussions. This work is supported by a Five-Year Fellowship from
the American Institute of Mathematics.

\section{The Khovanov bound}
\label{sec:KhBound}

This section contains the proof of Theorem~\ref{thm:strong} and some
immediate applications.

We first review a bit of Legendrian knot theory. A generic
Legendrian link in standard contact $\R^3$ projects in the $xz$
plane to a \textit{front}, whose only singularities are double
points and semicubical cusps, and which has no vertical tangencies.
A front can be resolved into a link diagram, which we call the
``desingularization'' of the front, by smoothing the cusps and
turning each double point into a crossing in which the strand of
more negative slope is the overstrand:
\raisebox{-0.05in}{\includegraphics[width=0.2in]{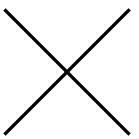}}
$\to$
\raisebox{-0.05in}{\includegraphics[width=0.2in]{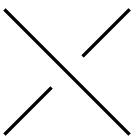}}.
For an oriented front $F$, let $w(F)$ denote the writhe of the
desingularization of $F$, and $c(F)$ half the number of cusps of
$F$; then the Thurston--Bennequin number of $F$ is given by $\tb(F)
= w(F)-c(F)$.

Via desingularization, any front $F$ has a Khovanov complex
$\CKh(F)$ whose homology is Khovanov homology. This complex is
bigraded, but in this section (except for the proof of
Corollary~\ref{cor:s-invt}) we will collapse one direction of the
grading and consider the single grading given by the quantum ($q$)
grading minus the homological ($t$) grading: $|\v| = p(\v) + w(F)$.
Here $\v$ is a tensor product of vectors $\v_\pm$, one associated to
each component of a resolution of the desingularization of $F$, and
$p$ is defined by setting $p(\v_\pm) = \pm 1$. With respect to this
grading, the homology $\HKh^*(F)$ of $\CKh^*(F)$ is an invariant of
the link type of the desingularized front.

We will be primarily concerned with a shifted version $\CKhshift$ of
the Khovanov complex, with the grading given by $|\v|' = p(\v)$.
Note that $\CKhshift(F)$ and its homology $\HKhshift(F)$ do not
depend on an orientation of $F$, and that $\HKh^*(F) =
\HKhshift^{*-w(F)}(F)$.

\begin{proposition}
$\HKhshift^*(F) = 0$ for $* < -c(F)$, and thus $\HKh^*(F) = 0$ for
$* < \tb(F)$.
\label{prop:vanishing}
\end{proposition}

Theorem~\ref{thm:strong} follows immediately from
Proposition~\ref{prop:vanishing}. We will prove
Proposition~\ref{prop:vanishing} by induction on the number of
double points of $F$. First we need to introduce a bit of
terminology.

We say that a front $F$ is $n$-\textit{vanishing} if $\HKhshift^*(F)
= 0$ for $* < n$; Proposition~\ref{prop:vanishing} states that any
front is $-c(F)$-vanishing. Given an unoriented front $F$ and a
double point $p$ in $F$, we can construct two new fronts
$\Res_0(F,p),\Res_1(F,p)$ which, in the desingularized picture, are
Khovanov's $0,1$-resolutions of $F$ at the crossing $p$; that is,
$\Res_0(F,p)$ and $\Res_1(F,p)$ are the fronts obtained from $F$ by
replacing the double point
\raisebox{-0.05in}{\includegraphics[width=0.2in]{figures/crossing.eps}}
at $p$ by
\raisebox{-0.05in}{\includegraphics[width=0.2in]{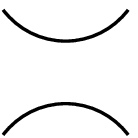}}
and
\raisebox{-0.05in}{\includegraphics[width=0.2in]{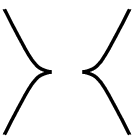}},
respectively. Note that $c(F) = c(\Res_0(F,p)) = c(\Res_1(F,p)) -
1$.

\begin{lemma}
There is a long exact sequence
\[
\xymatrix{ \HKhshift(\Res_0(F,p)) \ar[rr]^{(-1)} &&
\HKhshift(\Res_1(F,p))
\ar[dl] \\
& \HKhshift(F), \ar[ul] & }
\] where the top map lowers degree by $1$ and
the other maps preserve degree. \label{lem:exact}
\end{lemma}

\begin{proof}
This follows directly from the short exact sequence of complexes
\[
0\to \CKhshift(\Res_1(F,p)) \to \CKhshift(F) \to
\CKhshift(\Res_0(F,p))\to 0,
\]
which itself follows from the facts that $\CKhshift(F)$ without its
differential is the direct sum $\CKhshift(\Res_0(F,p)) \oplus
\CKhshift(\Res_1(F,p))$, and its differential preserves
$\CKhshift(\Res_1(F,p))$.
\end{proof}

\begin{lemma}
If $\Res_0(F,p)$ and $\Res_1(F,p)$ are $n$-vanishing for some $n$,
then so is $F$.
\label{lem:vanishing}
\end{lemma}

\begin{proof}
Immediate from the exactness of $\HKhshift^*(\Res_1(F,p)) \to
\HKhshift^*(F) \to \HKhshift^*(\Res_0(F,p))$ from
Lemma~\ref{lem:exact}.
\end{proof}

\begin{proof}[Proof of Proposition~\ref{prop:vanishing}]
We induct on the number of double points of $F$. If $F$ has no
double points, then it is an unlinked union of $n$ unknots, where $n
\leq c(F)$. An unknot has $\HKh^*$ supported in dimensions $-1$ and
$1$, and so $\HKh^*(F)$ is supported in dimensions between $-n$ and
$n$; since $w(F) = 0$, it follows that $F$ is $-n$-vanishing and
hence $-c(F)$-vanishing.

Now consider an arbitrary front $F$. Let $p$ be the double point of
$F$ farthest to the right. We consider two cases.

If the two strands emanating from the right of $p$ meet at a right
cusp (i.e., $p$ is part of a ``fish''), then the desingularization
of $F$ is the same as the desingularization of $\Res_0(F,p)$ except
for the addition of a negative kink. If we give $F$ any orientation
and $\Res_0(F,p)$ the induced orientation, then $w(F) =
w(\Res_0(F,p))-1$, and invariance under (topological) Reidemeister
move I implies that $\HKh^*(F) = \HKh^*(\Res_0(F,p))$ for all $*$.
Thus $\HKhshift^*(F) = \HKhshift^{*-1}(\Res_0(F,p))$ for all $*$,
and the induction assumption that $\Res_0(F,p)$ is $-c(F)$-vanishing
shows that $F$ is as well.

Otherwise, the two strands from the right of $p$ do not end at the
same cusp. By the induction assumption, $\Res_0(F,p)$ is
$-c(F)$-vanishing, and thus the induction step follows from
Lemma~\ref{lem:vanishing} if we can show that $\Res_1(F,p)$ is also
$-c(F)$-vanishing.

In terms of singularities, $\Res_1(F,p)$ replaces $p$ by two cusps,
a right cusp (opening to the left) and a left cusp (opening to the
right). Starting at the left cusp, follow the front $\Res_1(F,p)$ in
either direction until the front passes to the left of $p$; the
result is a zigzag path including the left cusp which does not cross
itself or the rest of the front. For consecutive cusps along this
zigzag, measure the difference between the $x$ coordinates of the
cusps, and let $c_1,c_2$ be the consecutive cusps for which this
difference is smallest. If we traverse the zigzag path in either
direction, $c_1$ and $c_2$ must be traversed both downwards or both
upwards by minimality. The zigzag between these two cusps comprises
a ``stabilization'' which can be eliminated to obtain another front
$F'$ which agrees with $\Res_1(F,p)$ outside of the stabilization;
see \figref{fig:cusp-elimination}.

\begin{figure}[ht!]\anchor{fig:cusp-elimination}
\centerline{
\small
\SetLabels
\E(.08*.32) $c_1$\\
\E(.25*.75) $c_2$\\
\endSetLabels
\AffixLabels{
\includegraphics[height=1in]{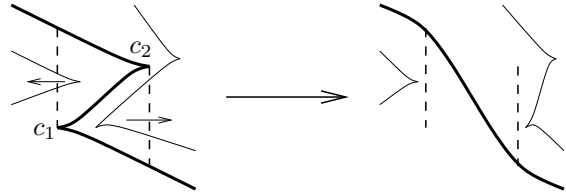}}}
\caption{
Eliminating cusps $c_1$ and $c_2$. One may need to first perform an
isotopy which pushes any part of the front with $x$ coordinates
between those of $c_1$ and $c_2$ out of the stabilization region.
}
\label{fig:cusp-elimination}
\end{figure}

Clearly $\Res_1(F,p)$ and $F'$ have the same $\HKhshift$ since they
have the same front desingularization. But now $c(F') =
c(\Res_1(F,p))-1 = c(F)$; by the induction hypothesis, $F'$ and thus
$\Res_1(F,p)$ is $-c(F)$-vanishing, as desired.
\end{proof}

The proof of Proposition~\ref{prop:vanishing} also yields a
sufficient condition for a front to maximize $\tb$ in its
topological class. Let the $0$-\textit{resolution} of a front be the
front which results from taking the $0$-resolution of each double
point:
\raisebox{-0.05in}{\includegraphics[width=0.2in]{figures/crossing.eps}}
$\to$
\raisebox{-0.05in}{\includegraphics[width=0.2in]{figures/crossing-res0.eps}}.
A $0$-resolution is \textit{admissible} if: (1) each component of
the $0$-resolution contains exactly two cusps; and (2) at each
resolved double point, the two arcs of the resolution belong to
different components of the $0$-resolution.

\begin{proposition}
Any front with admissible $0$-resolution maximizes
Thurston--Bennequin number in its topological link class; both weak
and strong Khovanov bounds are sharp in this case.
\label{prop:admissible}
\end{proposition}

\begin{proof}
Let $F$ be a front with admissible $0$-resolution. We wish to show
that $\HKhshift^{-c(F)}(F)$ has positive free rank. Label the double
points of $F$ from left to right $p_1,\dots,p_n$, and let
$F_0,F_1,\dots,F_n=F$ be the ``partial $0$-resolutions'' obtained
from $F$ such that $F_i$ is $F$ but with $0$-resolutions in place of
$p_{i+1},\dots,p_n$.

We prove by induction that $\HKhshift^{-c(F)}(F_i)$ has positive
free rank for all $i$. For $i=0$, $F_0$ is the $0$-resolution of
$F$, and by admissibility consists of $c(F)$ disjoint unknots; hence
$\HKhshift^{-c(F)}(F_0)=\Z$.

Now suppose that $\HKhshift^{-c(F)}(F_{i-1})$ has positive free
rank. By Lemma~\ref{lem:exact}, we have an exact sequence
\[
\HKhshift^{-c(F)}(F_i) \to \HKhshift^{-c(F)}(F_{i-1}) \to
\HKhshift^{-c(F)-1}(\Res_1(F_i,p_i)).
\]
By admissibility, the two strands emanating to the right of $p_i$ in
$F_i$ do not meet at a cusp, and so the proof of
Proposition~\ref{prop:vanishing} implies that $\Res_1(F_i,p_i)$ is
$-c(F)$-vanishing. Thus $\HKhshift^{-c(F)}(F_i) \to
\HKhshift^{-c(F)}(F_{i-1})$ is a surjection, and the induction step
is complete.
\end{proof}

Proposition~\ref{prop:admissible} immediately implies that the
Khovanov bound is sharp for closures of positive braids, i.e.,
braids consisting of positive products of elementary braid
generators, with the convention that the right handed trefoil is the
closure of a positive braid. Note that $\overline{\tb}$ was already
known in this case, as Bennequin's original bound is sharp.

We remark that an admissible $0$-resolution of a front constitutes
an ungraded normal ruling (or proper decomposition) of the front in
the sense of \cite{bib:ChP,bib:Fu}. It follows that the
Chekanov--Eliashberg differential graded algebra \cite{bib:ChP} for
the front has an augmentation, and in particular that the front is
not Legendrian isotopic to a stabilization (i.e., a front with a
small zigzag). Proposition~\ref{prop:admissible} strengthens the
result that the front is not destabilizable.
%

Note that since the original release of this paper, Rutherford
\cite{bib:Ru} has obtained a significant generalization of
Proposition~\ref{prop:admissible}: any front admitting an ungraded
normal ruling maximizes Thurston--Bennequin number in its
topological link class.

We will use Proposition~\ref{prop:admissible} in
Section~\ref{sec:alternating} to calculate $\overline{\tb}$ for
alternating links by showing that any alternating link has a front
with admissible $0$-resolution. It is not known what class of knots
has a front with admissible $0$-resolution; alternating links and
closures of positive braids fall into this category, while links for
which the Khovanov bound is not sharp (e.g., many negative torus
knots) do not.

To conclude this section, we use the weak Khovanov bound to deduce
the $s$-invariant bound on $\overline{\tb}$.

\begin{proof}[Proof of Corollary~\ref{cor:s-invt}]
Recall the definition of the $s$-invariant from \cite{bib:Ras}:
there is a spectral sequence from the bigraded Khovanov homology
$\HKh(K) \otimes \Q$ to another knot homology $\HKh'(K)$ introduced
by Lee \cite{bib:Lee}, and $\HKh'(K)=\Q\oplus\Q$ has a summand in
each of quantum gradings $s(K) \pm 1$. Lee gives an explicit set of
generators of $\HKh'(K)$, called ``canonical generators'' in
\cite{bib:Ras}; it is easy to see that both canonical generators
have homological grading $0$. It follows that $\HKh'(K)$ and thus
$\HKh^{*,*}(K)\otimes\Q$ is nonzero in bidegree $(s(K)\pm 1,0)$. Now
apply the weak Khovanov bound.
\end{proof}


\section{Alternating links}
\label{sec:alternating}

In this section, we show that the Khovanov bound is sharp for
alternating links. It is well known that there is a correspondence
between alternating links and planar graphs. For our purposes, a
planar graph is a graph embedded in the plane which may have more
than one edge connecting a pair of vertices, but has no edges
connecting a vertex to itself. We consider planar graphs up to
isotopies of the plane.

\begin{definition}
A \textit{reduced planar graph} is a $1$-connected planar graph;
that is, it is connected and cannot be disconnected by removing one
edge. In particular, it has no vertices of valence $1$.
\end{definition}

There is a standard way to obtain a reduced planar graph from a
reduced alternating link diagram. (Recall that an alternating
diagram is reduced if there is no crossing whose removal splits the
diagram into two disjoint parts; any alternating link has a reduced
alternating diagram.) Such a diagram divides the plane into a number
of components; color these components in a checkerboard fashion so
that near any crossing
\raisebox{-0.05in}{\includegraphics[width=0.2in]{figures/crossing-alt.eps}},
the coloring takes the form
\raisebox{-0.06in}{\includegraphics[width=0.2in]{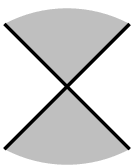}}.
The vertices of the planar graph correspond to black regions, and
the edges are diagram crossings where two black regions meet.

We will need another way to depict a planar graph.

\begin{definition}
A \textit{Mondrian diagram} consists of a set of disjoint horizontal
line segments in the plane, along with a set of disjoint vertical
line segments, each of which begins and ends on a horizontal segment
and does not intersect any other horizontal segment.
\end{definition}

Any Mondrian diagram yields a planar graph by contracting the
horizontal segments to points.

\begin{proposition}
All reduced planar graphs are the contraction of some Mondrian
diagram. \label{prop:Mondrian}
\end{proposition}

In fact, Proposition~\ref{prop:Mondrian} holds for arbitrary planar
graphs; this generalization is easy to establish from
Proposition~\ref{prop:Mondrian}, but we need only the reduced case
for our purposes.

We delay the proof of Proposition~\ref{prop:Mondrian} until the end
of this section. First we apply it to the maximal
Thurston--Bennequin number for alternating links.

\begin{figure}[ht!]\anchor{fig:63-front}
\centerline{
\includegraphics[width=4in]{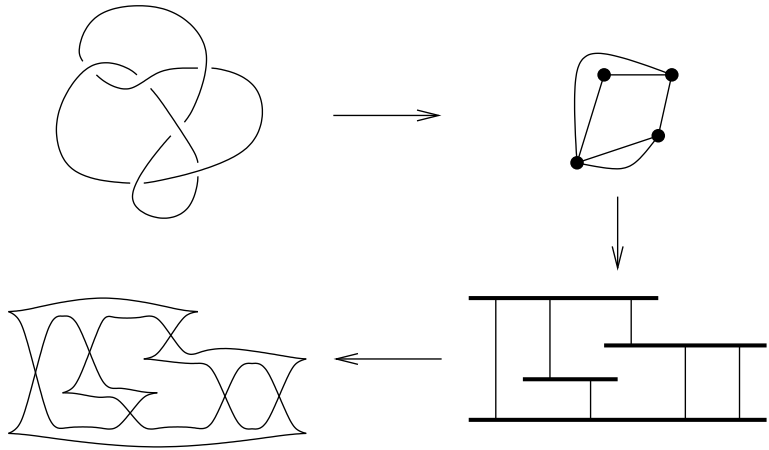}
}
\caption{
Using a Mondrian diagram to construct a Legendrian form for the knot
$6_3$. Clockwise, from left: a reduced alternating diagram for
$6_3$; the corresponding planar graph; a Mondrian diagram
contracting to this graph; a front for $6_3$.
}
\label{fig:63-front}
\end{figure}

Let $K$ be an alternating link; this has a reduced alternating
diagram which gives rise to a planar graph. Consider a Mondrian
diagram whose contraction is this graph. If necessary, extend the
ends of the horizontal segments of this diagram slightly so that no
vertical segment ends at an endpoint of a horizontal segment. We can
now turn this Mondrian diagram into a front as follows: replace each
horizontal segment by a ``pair of lips'' front for the unknot, and
delete from these fronts a neighborhood of each intersection with a
vertical segment; then replace each vertical segment by an X. This
front represents a Legendrian link of the topological type of $K$.
See \figref{fig:63-front} for an example.

We are now in a position to prove the sharpness of the Khovanov
bound for alternating links.

\begin{proof}[Proof of Theorem~\ref{thm:alternating}]
The fact that the weak (and hence also the strong) Khovanov bound is
sharp for alternating knots follows from
Proposition~\ref{prop:admissible} and the observation that the front
constructed above from a Mondrian diagram is admissible. To complete
the proof of the proposition, we need to establish that $\mindeg_q
\Kh_K(q,t/q) = \mindeg_q V_K(q) + \sigma(K)/2-1$ for alternating
nonsplit links $K$. This follows easily from the identity
$\Kh_K(q,-1) = (q+q^{-1})V_K(q^2)$, along with a result of Lee
\cite{bib:Lee}:
\[
\Kh_K(q,t) = q^{\sigma(K)} (q^{-1}+q+(q^{-1}+tq^3)\Kh'_K(tq^2)),
\]
where $\Kh'_K$ is some Laurent polynomial. Note that our sign
convention for $\sigma$ is the opposite of Lee's.
\end{proof}

We remark that, by work of Murasugi or Thistlethwaite, the
expression from Theorem~\ref{thm:alternating} for $\overline{tb}(K)$
when $K$ is alternating and nonsplit can be rewritten as follows:
\[
\overline{tb}(K) = -c_-(K) + \sigma(K) - 1,
\]
where $c_-(K)$ is the number of negative crossings in any reduced
alternating diagram for $K$. As another side note,
Theorem~\ref{thm:alternating}, together with the fact that the
Kauffman bound is also sharp for alternating links \cite{bib:Ru},
implies that for $K$ alternating,
\[
\mindeg_q \Kh_K(q,t/q) = \mindeg_a F_K(a,z) - 1,
\]
where $F_K(a,z)$ is the Kauffman polynomial of $K$, normalized so
that the Kauffman polynomial of the unknot is $1$.

\begin{figure}[ht!]\anchor{fig:mondrian-ex}
\centerline{
\includegraphics[height=1.5in]{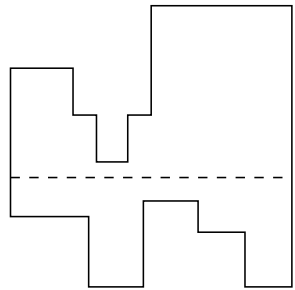}
}
\caption{
A strong region, with spine given by the dashed line. For an example
of a strong Mondrian diagram, see Figure~\ref{fig:63-front}.
}
\label{fig:mondrian-ex}
\end{figure}

The rest of this section is devoted to a proof of
Proposition~\ref{prop:Mondrian}. We will actually prove a slightly
stronger statement. The complement of a Mondrian diagram in the
plane consists of several connected components, one unbounded and
the rest bounded. Call a bounded component \textit{strong} if the
interior $R$ of its closure has the following property: there is a
horizontal line segment $L$ in $R$ such that every nonempty vertical
slice of $R$ (the intersection of $R$ with a vertical line) is
connected and intersects $L$. In this case, $L$ is called a
\textit{spine} of the region. See \figref{fig:mondrian-ex}. A
Mondrian diagram is called strong if each bounded component of its
complement is strong.

We will prove that any reduced planar graph is the contraction of a
strong Mondrian diagram. Let $G$ be a reduced planar graph; we
construct a strong Mondrian diagram for $G$ from the outside in.

Define the \textit{boundary} of a planar graph to be the subgraph
consisting of vertices and edges abutting the unbounded region in
the complement of the graph. Call a planar graph an \textit{enhanced
cycle} if its boundary is a cycle (with no repeated vertices); an
enhanced cycle consists of this cycle, along with some number of
edges inside the cycle which connect vertices of the cycle.

We will construct strong Mondrian diagrams for a sequence of
subgraphs of $G$ which build up to $G$.

\begin{figure}[ht!]\anchor{fig:cycle}
\centerline{
\includegraphics[height=2.4in]{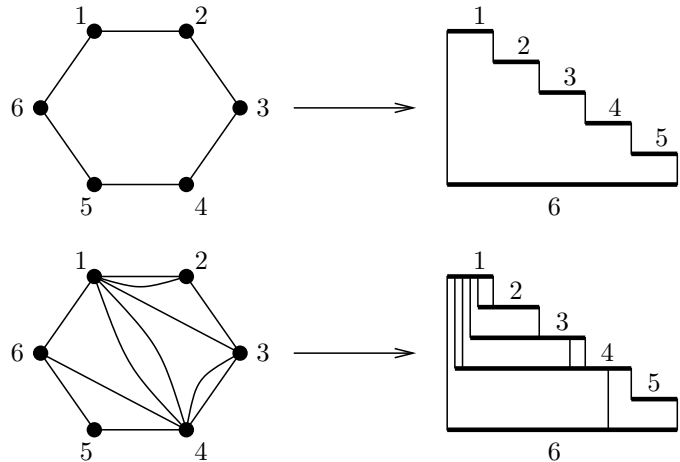}
}
\caption{
Step-shaped Mondrian diagram for a cycle (top), and resulting strong
Mondrian diagram for an enhanced cycle (bottom)}
\label{fig:cycle}
\end{figure}

\begin{step}
Any planar cycle. Here a strong Mondrian diagram is given by a
step-shaped construction as shown in \figref{fig:cycle}.
\end{step}

\begin{step}
An enhanced cycle $C$. Number the vertices of $C$ $1,\dots,n$ in
clockwise order around the boundary cycle $\tilde{C}$. From the
previous step, $\tilde{C}$ is the contraction of a step-shaped
Mondrian diagram whose ``steps'' from top to bottom correspond in
order to $1,\dots,n$. There is now an essentially unique way to
expand this Mondrian diagram to give one for $C$; see
\figref{fig:cycle}. More precisely, for an edge in $C\setminus
\tilde{C}$ joining vertices $i$ and $j$ with $i<j$, drop a
perpendicular from the step $i$ to the level of $j$. Do this for
each edge in $C\setminus \tilde{C}$, and arrange the perpendiculars
so that the lengths of the perpendiculars dropped from any
particular step $i$ are nonincreasing as we view the perpendiculars
from left to right. Now extend each step leftwards just as far as it
needs to go to meet each perpendicular which ends at its height. The
planarity of $\tilde{C}$ implies that the resulting diagram is a
Mondrian diagram; its strongness follows from the easily checked
fact that each bounded region in its complement is step-shaped.
\end{step}

For the next steps, label the vertices of $G$ $v_1,\dots,v_m$ in
such a way that $v_1,\dots,v_k$ for some $k$ are the vertices on the
boundary of $G$, and the subgraph of $G$ induced by $v_1,\dots,v_l$
is connected for all $l\geq k$.

\begin{figure}[ht!]\anchor{fig:tree}
\centerline{
\includegraphics[width=4.5in]{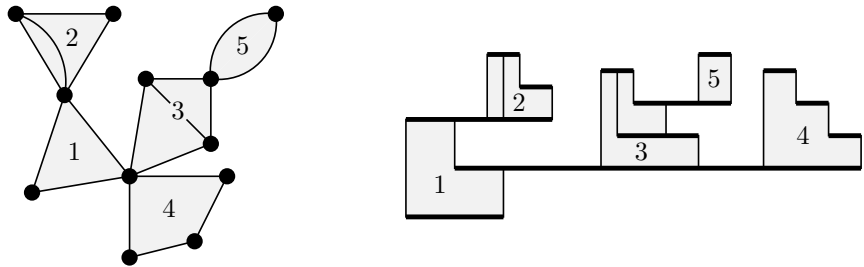}
}
\caption{
Joining Mondrian diagrams for enhanced cycles (labeled $1,2,3,4,5$)
}
\label{fig:tree}
\end{figure}

\begin{step}
The subgraph $G_k$ of $G$ induced by $v_1,\dots,v_k$. Since the
boundary of $G$ consists of an edge-disjoint union of cycles, $G_k$
consists of an edge-disjoint union of enhanced cycles, joined in
treelike fashion by their common vertices. We can build a strong
Mondrian diagram for $G_k$ by starting with the strong Mondrian
diagram for one of the enhanced cycles, extending rightward each
horizontal segment corresponding to a vertex to which other enhanced
cycles are attached, placing strong Mondrian diagrams for these
attached enhanced cycles with base on the extended segments, and
continuing in this fashion until all of $G_k$ has been constructed.
Here we use the fact that any vertex of an enhanced cycle can serve
as the base of its step-shaped Mondrian diagram. See
\figref{fig:tree} for a (probably clearer) pictorial example.
\end{step}

\begin{step}
\label{step:last}
The subgraph $G_l$ of $G$ induced by $v_1,\dots,v_l$, for $k<l\leq
m$. This proceeds by induction on $l$. Assume that we have a strong
Mondrian diagram which contracts to $G_{l-1}$; we construct one for
$G_l$. Note that $G_l$ is $G_{l-1}$ along with vertex $v_l$ and the
edges from $v_l$ to $G_{l-1}$; we may assume that there are at least
two such edges, since if there is only one, we can double it and
then remove the double at the end of the induction process.

Let $R$ be the component of the complement of $G_{l-1}$ which
contains $v_l$, corresponding to a strong region $\tilde{R}$ in the
strong Mondrian diagram for $G_{l-1}$. Fix a spine of $\tilde{R}$.
For each edge from $v_l$ to a vertex on the boundary of $R$, locate
the (unique) horizontal boundary component of $\tilde{R}$
corresponding to this vertex, and drop a vertical perpendicular from
this horizontal segment to the spine. We may assume that no two
perpendiculars land on the same point of the spine. The subset of
the spine beginning at the foot of the leftmost perpendicular and
ending at the foot of the rightmost perpendicular will contract to
the vertex $v_l$. See \figref{fig:mondrian-induction}.

\begin{figure}[ht!]\anchor{fig:mondrian-induction}
\centerline{
\includegraphics[width=5in]{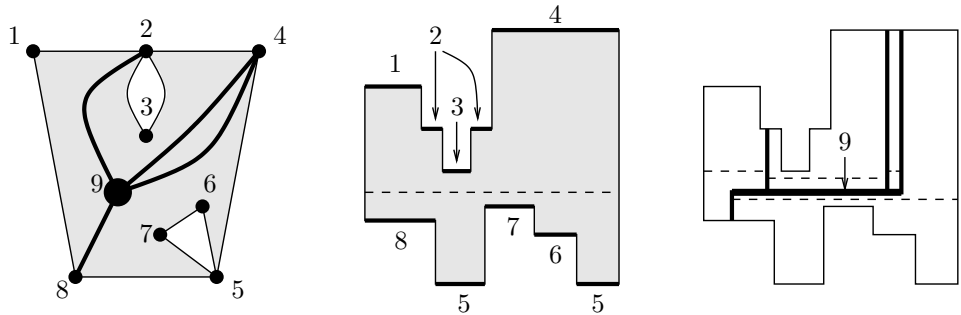}
}
\caption{
Induction step in the construction of a strong Mondrian diagram. We
wish to add vertex $9$ in the region $R$ in the planar graph (left)
to the region $\tilde{R}$ in the Mondrian diagram (center). A
portion of the spine (dashed line, center) becomes the horizontal
segment corresponding to vertex $9$. The result is the figure on the
right (with spines for each region given by the dashed lines). }
\label{fig:mondrian-induction}
\end{figure}

The Mondrian diagram for $G_l$ now consists of the Mondrian diagram
for $G_{l-1}$, along with this shortened spine and the
perpendiculars to it. To see that it is strong, note that for each
region into which $\tilde{R}$ is split, a horizontal slice of the
region either just above or just below the spine of $\tilde{R}$
constitutes a spine for the region.
\end{step}

At the end of the induction for Step~\ref{step:last}, we have a
strong Mondrian diagram for $G$, possibly with some doubled edges.
By removing the extraneous vertical segments corresponding to these
doubled edges, we obtain a strong Mondrian diagram for $G$, as
desired.

\section{Computations}
\label{sec:computations}

We now discuss various computations of the Khovanov bound, its
sharpness, and its relation to other known bounds on
$\overline{\tb}$. We use the standard knot notation as shown in
Rolfsen for knots with at most $10$ crossings, and
Dowker--Thistlethwaite notation for knots with $11$ crossings, and
denote mirroring by an overline; for instance, $11n_{24}$ denotes
the $24^\textrm{th}$ nonalternating $11$-crossing knot in the
Dowker--Thistlethwaite enumeration, and $\overline{11n_{24}}$ is its
mirror. Many computations in this section were assisted by the
program \textit{Knotscape} \cite{bib:Knotscape}, the
\textit{Mathematica} package \texttt{KnotTheory\`{}} \cite{bib:BN},
the online Table of Knot Invariants \cite{bib:KnotInfo}, and the
Khovanov homology data of Shumakovitch \cite{bib:ShuData}.

The Khovanov bound (either strong or weak) is quite effective in
calculating maximal Thurston--Bennequin number. It is sharp for all
knots with $9$ or fewer crossings, and in particular resolves the
one unknown value in the table of $\overline{\tb}$ from
\cite{bib:NgTB}.

\begin{proposition}
$\overline{\tb}(\overline{9_{42}}) = -5.$
\end{proposition}

For $10$-crossing knots, Theorem~\ref{thm:alternating} allows us to
restrict our attention to the $42$ prime nonalternating knots and
their mirrors. Hand-drawn front diagrams by the author (see the
author's web page for details) yield the following result.

\begin{proposition}
The Khovanov bound (either strong or weak) is sharp for all knots
with $10$ or fewer crossings, with the exception of
$\overline{10_{124}}$, for which
$\overline{\tb}(\overline{10_{124}}) = -15$ while the Khovanov bound
is $-14$, and the possible exception of $\overline{10_{132}}$, for
which we have $-1 \leq \overline{\tb}(\overline{10_{132}}) \leq 0$.
\label{prop:10crossing}
\end{proposition}

\begin{figure}[ht!]\anchor{fig:10132-front}
\centerline{
\includegraphics[height=1.35in]{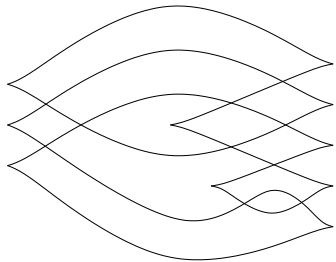}
}
\caption{
A possible maximal-$\tb$ representative for $\overline{10_{132}}$?
}
\label{fig:10132-front}
\end{figure}

\noindent For the negative torus knot $\overline{10_{124}}=T(5,-3)$,
the computation of $\overline{\tb}$ comes from \cite{bib:EH}. A
diagram of $\overline{10_{132}}$ with $\tb=-1$ is given in
\figref{fig:10132-front}, while the Khovanov bound is $0$. A
table of values for $\overline{\tb}$ for prime knots up to $10$
crossings, deduced from Proposition~\ref{prop:10crossing}, is
available online at \cite{bib:KnotInfo}.

Negative torus knots provide a number of further examples of
nonsharpness of the Khovanov bound. The Khovanov bound (strong or
weak) is sharp for negative torus knots $T(2n+1,-2)$, $T(4,-3)$, and
$T(7,-5)$, but for no other examples where the author has computed
it. We note that the Kauffman bound is sharp for $T(p,-q)$ when
$p>q$ and $q$ is even \cite{bib:EF}; when $p>q$ and $q$ is odd, it
appears that the Khovanov bound gives a better estimate than
Kauffman, though neither is sharp in general.

As a particular example, consider $T(5,-4)$. Here the strong and
weak Khovanov bounds for $\overline{\tb}(T(5,-4))$ are
different,\footnote{This seems to be a relatively unusual situation.
There are no knots with crossing number $13$ or less for which the
strong and weak Khovanov bounds disagree.} given by $-19$ and $-18$,
respectively \cite{bib:Shu}; the difference is the torsion group
$\HKh^{-29,-10}(T(5,-4)) = \Z/2$. Here neither Khovanov bound is
sharp, while the Kauffman bound does give the sharp value
$\overline{\tb} = -20$.

We see from this example that the Khovanov and Kauffman bounds are
incommensurate. However, for small knots at least, the Khovanov
bound often seems to be better when the two disagree. There are $19$
knots with $11$ or fewer crossings for which the bounds disagree:
$\overline{8_{19}}$, $\overline{9_{42}}$, $\overline{10_{124}}$,
$\overline{10_{128}}$, $\overline{10_{136}}$, $11n_{20}$,
$\overline{11n_{24}}$, $\overline{11n_{27}}$, $\overline{11n_{37}}$,
$\overline{11n_{50}}$, $\overline{11n_{61}}$, $\overline{11n_{70}}$,
$\overline{11n_{79}}$, $\overline{11n_{81}}$, $11n_{86}$,
$\overline{11n_{107}}$, $\overline{11n_{126}}$,
$\overline{11n_{133}}$, and $\overline{11n_{138}}$. For all of
these, the Khovanov bound is better. Of the $46$ $12$-crossing knots
where the bounds disagree, Kauffman is better for one ($12n_{475}$,
also the only knot with $12$ or fewer crossings for which the
HOMFLYPT bound is stronger than Khovanov), and Khovanov is better
for the rest.



\Addresses\recd
\end{document}